\documentclass[twocolumn]{autart}    
\usepackage{natbib}
\usepackage{amsmath}
\usepackage{amssymb}
\usepackage{glossaries}
\usepackage[T1]{fontenc}
\newacronym{lti}{LTI}{linear time-invariant}
\newacronym{mimo}{MIMO}{multiple-input, multiple-output}
\newacronym{pe}{PE}{persistency of excitation}
\usepackage{graphicx}          

\begin{document}

\begin{frontmatter}
\title{Data-driven representation of continuous-time linear time-invariant behaviors via Hankel integral operators}

\author[Delft]{Fritz A. Engeln}\ead{f.a.engeln@tudelft.nl}, 
\author[Delft]{Marta Zag{\'o}rowska}\ead{m.a.zagorowska@tudelft.nl},
\author[Delft]{Jan-Willem van Wingerden}\ead{j.w.vanwingerden@tudelft.nl}

\address[Delft]{Delft Center for Systems and Control, Delft University of Technology, Mekelweg 2, 2628CD Delft, The Netherlands}
          

\begin{abstract}                          
    We present a data-driven representation of the finite-horizon behavior of deterministic controllable continuous-time LTI systems. We first introduce a notion of persistency of excitation for continuous-time signals that does not require the signal to be differentiable. We then show that, when the input is persistently exciting in this sense, the closure of the range of a Hankel integral operator built from a single input-output trajectory equals the behavior of the system over a finite time horizon. This result provides a nonparametric representation of system trajectories directly in terms of measured data, without requiring an explicit parametric model of the system.
\end{abstract}

\end{frontmatter}

\section{Introduction}
\label{sec:Introduction}
A central challenge in control engineering is that an accurate parametric model of the system is often not available. System identification addresses this by estimating such a model from measured input-output data \citep{ljung1999system}. Data-driven control approaches offer an alternative by exploiting measured data directly to design controllers without explicitly estimating a parametric model \citep{8795639,FAVOREEL19994004}. However, a question arises whether a given data trajectory contains enough information to represent the behavior of the underlying system. The notion of \gls{pe} addresses this question. Like much of data-driven control, it has been developed in discrete time, where it is expressed as a rank condition on Hankel matrices constructed from shifted data trajectories. The result by \cite{WILLEMS2005325}, known as the \textit{fundamental lemma}, states that if the input of a controllable discrete-time \gls{lti} system is \gls{pe} of sufficient order, every finite-length input-output trajectory resides in the column span of a matrix constructed from the input-output data.

Recent developments in areas such as event-triggered and networked control have renewed interest in continuous-time control methods, which allow the evolution of the system and the implementation of control actions to be analyzed without relying on periodic sampling \citep{HERTNECK2026101387}. This motivates the development of theoretical tools for data-driven control and system identification of continuous-time \gls{lti} systems. However, the notion of \gls{pe} does not extend directly to continuous time. In discrete time, the finite-length behavior of an \gls{lti} system is a finite-dimensional linear subspace, which can be spanned by the finitely many columns of a Hankel matrix. In continuous time, signals are functions rather than finite sequences, and the space of trajectories on a finite interval is infinite-dimensional. A finite collection of data trajectories therefore cannot span it, and finite-dimensional data matrices must be replaced by operators on function spaces.

The informativity of data for continuous-time systems was first studied in the context of adaptive control and identification. \cite{MAREELS1984159} and \cite{SHIMKIN1987225} define conditions on the input data that guarantee convergence of parameter estimates for continuous-time \gls{lti} systems. In both works, the condition is built from the derivatives of the input. These conditions are sufficient for system identification, but they do not characterize the subspace in which the input and output functions of the system reside.\\
A continuous-time counterpart of the \textit{fundamental lemma} was first given by \cite{9992347}. They split the data into segments of equal length and, at each time instant, arrange the samples taken one segment length apart in a time-dependent Hankel matrix. The data are called \gls{pe} if this matrix has full row rank at every time instant. The condition is imposed on the input and state data jointly, so it depends on the system and its initial state. The coefficients that combine the shifted data are not constant. They solve a differential equation that requires measured state data and is well defined only if the generated input is differentiable. \cite{10540062} remove the need for state measurements by working with the input and output together with their derivatives. The coefficients again solve a differential equation, the derivatives of the data up to the system lag must be available, and only trajectories with a continuously differentiable input are described. A related approach is presented by \cite{10886882}.
\\
A different notion of \gls{pe} is introduced by \cite{9882335}. An input is \gls{pe} if its derivatives up to a given order are linearly independent. This is closely related to the condition of \cite{MAREELS1984159}, and implies that the input does not satisfy any linear differential equation with constant coefficients of given order. For controllable systems, such inputs allow one to recover the differential equations of the system from a single input-output trajectory on a finite interval. This result concerns the identification of the system and does not describe its trajectories in terms of the data. The input has to be differentiable, and recovering the equations requires the derivatives of both the input and the output. Building on this notion, \cite{RAPISARDA2023105603} represent the trajectories of the system by the coefficients of a Chebyshev series, and \cite{SCHMITZ2024105950} describe them, together with their derivatives, by a matrix built from integrals of the data and their derivatives and a time-dependent coefficient function. The Chebyshev-based representation has been used for data-driven iterative learning control \citep{10384133} and predictive control \citep{10886618}. The recent work by \cite{WOLSKI2026106594} weights the data with one modulating function at a finite number of time shifts. The derivatives are moved onto the test function, so they do not have to be measured. The resulting matrix of samples characterizes the smooth trajectories of a controllable system. However, the excitation condition is the one of \cite{9882335}, which requires differentiable inputs.\\
A sufficient \gls{pe} condition for identifying a system differs from a condition under which recorded data span the finite-length behavior of the system. In discrete time, the \textit{fundamental lemma} of \cite{WILLEMS2005325} gives a condition under which finite linear combinations of shifted data trajectories span the space of all trajectories of a given length. This condition is generally stronger than what identifiability requires. One could instead identify the system and then generate trajectories from the identified model to span its behavior. In continuous time, \cite{9882335} give a condition for identifiability and then span the behavior using polynomial expansions of functions that belong to the system \citep{RAPISARDA2023105603}. In contrast, we introduce a notion of \gls{pe} under which the behavior is spanned directly by integrals over shifted windows of a single measured \gls{pe} input.

The existing works rely on excitation conditions that require (piecewise) differentiable inputs. A characterization of the finite-horizon behavior, under a data informativity condition that requires no differentiability, has so far been missing. The contributions of this paper are:
\begin{enumerate}
    \item[(i)] We introduce a notion of \gls{pe} for \gls{mimo} continuous-time \gls{lti} systems that is defined on the input data only. It requires that the cross-correlation of the input with a test function and the input itself cannot be reproduced by an \gls{lti} system of the same order driven by the same input. Unlike the existing conditions, it requires no differentiability of the input.
    
    \item[(ii)] We show that, for a controllable system and an input that is \gls{pe}, the closure of the range of a Hankel integral operator built from shifted input-output data equals the behavior of the system over a finite horizon. Every element of the range is a system trajectory, and the Hankel operator requires neither state measurements nor derivatives of the data. The behavior is thus represented directly in terms of the data, as a closed subspace of the Hilbert space of square-integrable functions.
\end{enumerate}
The remainder of this paper is organized as follows. Section~\ref{sec:Preliminaries} introduces the preliminaries. In Section~\ref{sec:main} we introduce the notion of \gls{pe} for continuous-time systems and present our main result and its proof. Section~\ref{sec:Example} presents an illustrative example and Section~\ref{sec:Conclusion} concludes the paper.

\textit{Notation:} Let $L^2(\mathcal I;\mathbb R^n)$ denote the Hilbert space of square-integrable functions from the interval $\mathcal I$ to $\mathbb R^n$, equipped with the standard inner product. For a normed space $X$ with dual space $X^\ast$, we write $\langle x,x^\ast\rangle_{X,X^\ast}=x^\ast(x)$ for the canonical pairing of $x\in X$ and $x^\ast\in X^\ast$. If $X=H$ is a Hilbert space, the canonical pairing coincides with the inner product on $H$. Let $\delta$ denote the Dirac delta at $0$. Function composition is denoted by $\circ$, e.g., $(f\circ g)(\cdot) = f(g(\cdot))$. The indicator function of a set $\mathcal{S}$ is denoted by $\mathbf{1}_{\mathcal{S}}$.
\section{Preliminaries}
\label{sec:Preliminaries}
We consider the deterministic continuous-time \gls{lti} system
\begin{equation}
\label{eqn:LTIsys}
    \begin{aligned}
        \dot{x}(t) &= Ax(t) + Bu(t)\\
        y(t) &= C x(t) + Du(t),
\end{aligned}
\end{equation}
where $A\in \mathbb{R}^{n \times n}$, $B\in \mathbb{R}^{n\times n_u}$, $C\in \mathbb{R}^{n_y \times n}$, and $D\in\mathbb{R}^{n_y\times n_u}$ are the system, input, output, and direct feedthrough matrix, respectively. The vectors $x(t) \in \mathbb{R}^{n}$, $u(t)\in\mathbb{R}^{n_u}$ and $y(t)\in\mathbb{R}^{n_y}$ are the state, input, and output at time $t$, respectively.

\begin{minipage}{\linewidth}
    \begin{assum}
        Assume that the input, output, and state functions of the system in \eqref{eqn:LTIsys} are square integrable on finite intervals $\mathcal{I}\subset \mathbb{R}$, i.e.,
        \begin{align*}
            u \in L^2(\mathcal{I};\mathbb{R}^{n_u}),\quad y \in L^2(\mathcal{I};\mathbb{R}^{n_y}),\quad x \in L^2(\mathcal{I};\mathbb{R}^{n}).
        \end{align*}
    \end{assum}
\end{minipage}
\begin{assum}
\label{ass:control}
    The system in \eqref{eqn:LTIsys} is controllable, i.e., the controllability matrix
    \begin{align}
    \label{eqn:controllability}
        \mathcal{C} &= \begin{bmatrix}B&AB&\cdots & A^{n -1}B\end{bmatrix}
    \end{align}
    has full rank, such that $\operatorname{rank}(\mathcal{C})=n$.    
\end{assum}
For a given horizon length $\ell>0$, we define the length-$\ell$ input-output behavior as
\begin{equation}
\label{eqn:IOBehavior}
    \begin{aligned}
    \mathcal{B}_\ell &= \Big\{ (u,y):\exists x(0) \mathrm{\, s.t. \, \eqref{eqn:LTIsys}\,holds\,for\, a.e.\,} t\in [0,\ell]\Big\}\\
    &\subseteq L^2([0,\ell];\mathbb{R}^{n_u+n_y}).
\end{aligned}
\end{equation}
Analogously, we define the input-state behavior as
\begin{equation}
    \begin{aligned}
    \mathcal{B}_\ell^x &= \Big\{ (u,x): \mathrm{\;\eqref{eqn:LTIsys} \;holds \;for\; a.e.\;} t\in [0,\ell]\Big\}\\
    &\subseteq L^2([0, \ell];\mathbb{R}^{n_u+n}).
\end{aligned}
\end{equation}
Throughout the paper, we consider the input, output, and state functions $u_T$, $y_T$, and $x_T$ of length $T$ generated by the system in \eqref{eqn:LTIsys} such that
\begin{align}
    (u_T, y_T) \in \mathcal{B}_T \qquad \text{and} \qquad  (u_T, x_T) \in \mathcal{B}_T^x.
    \label{eqn:data}
\end{align}
Only the input-output functions $(u_T, y_T)$ are assumed to be measured. The state trajectory $x_T$ is used solely as a tool in the proofs and is not required to construct the data-driven representation. If the pair $(A,C)$ is not observable, $x_T$ is not uniquely determined by $(u_T, y_T)$; in that case, $x_T$ denotes any state trajectory that, together with $u_T$ and $y_T$, satisfies \eqref{eqn:LTIsys}.
For some $\ell$ such that $0<\ell<T$, we introduce the family of length-$\ell$ translates of $u_T$
\begin{align*}
    \mathcal{F}_{u_T,\ell}&=\left\{u\! :\! \exists s \! \in \![0,T-\ell] \,\text{s.t.} \,u(t)\!=\!u_T(t\!+\!s)\ \forall t\in[0,\ell]\right\}
    \\ &\subset L^2([0,\ell];\mathbb{R}^{n_u})
\end{align*}
whose elements are the restrictions of $u_T$ to windows of length $\ell$ starting at time $s$. Note that the family of translates itself is not a vector space. It is simply a set of functions inside the Hilbert space. However, the linear span of the family of translates forms a subspace of $L^2$.

\section{Main result}
\label{sec:main}
To quantify how informative an input signal is for the $n^{th}$-order system in \eqref{eqn:LTIsys}, we introduce the following notion of \gls{pe} for continuous-time signals:
\begin{defn}[Persistency of excitation]
\label{def:PE}
    Let $0 < \ell < T$. A signal $u_T \in L^2([0,T];\mathbb{R}^{n_u})$ is said to be
    $\ell$-PE of order $n$ if, for every pair $(v,q) \in L^2([0,\ell];\mathbb{R}^{n_u}) \times \mathbb{R}^{n_u}$ with $(v,q) \neq (0,0)$, the function
    \begin{align}
    \label{eqn:integralPE}
      y_{v,q}(s) \! = \!\int_0^\ell \! v(t)^\top u_T(t+s)\mathrm{d}t + q^\top u_T(s), \; s \in [0,T-\ell],
    \end{align}
    is not the output of any $n^{th}$-order strictly proper LTI system driven by $u_T$ for any initial state.
\end{defn}
The function $y_{v,q}$ combines the cross-correlation of a test function $v$ with the sliding window $u_T(\cdot+s)$ and the instantaneous value $q^\top u_T(s)$. Definition~\ref{def:PE} requires that no such combination can be reproduced by an $n^{th}$-order system driven by the same input. The condition involves only the input data and the order $n$, not the system matrices, and it requires no differentiability of $u_T$. An example of such a signal is given in Section~\ref{sec:Example}.

The following lemma shows that \gls{pe} implies that the translates of the
input are rich enough to approximate any function on the window $[0,\ell]$.
\begin{lem}
\label{lem:denseinput}
    Let $u_T\in L^2([0,T];\mathbb{R}^{n_u})$ be $\ell$-PE of any order. Then, the family of its length-$\ell$ translates $\mathcal{F}_{u_T,\ell}$ is complete\footnote{A family in a Hilbert space
    is complete if its linear span is dense in that space \citep{RudinRealAndComplex}.} in $L^2([0,\ell];\mathbb{R}^{n_u})$, i.e.,
    \begin{align}
        \overline{\operatorname{span}}(\mathcal{F}_{u_T,\ell}) &= L^2([0,\ell];\mathbb{R}^{n_u}).
    \end{align}
\end{lem}
\begin{pf}
    Suppose, for contradiction, that the family is not complete. Then its closed span is a proper closed subspace of $L^2([0,\ell];\mathbb{R}^{n_u})$, so there exists a
    nonzero $v$ orthogonal to every translate, i.e.,
    \begin{equation}
        \begin{aligned}
        y_{v,0}(s) &= \int_0^\ell v(t)^\top u_T(t+s)\mathrm{d}t \\
        &= \langle u_T(t+s), v(t)\rangle = 0
        \qquad \forall s \in [0,T-\ell].
    \end{aligned}
    \end{equation}
    The zero function is the output of a strictly proper LTI system. This contradicts Definition~\ref{def:PE} with $q = 0$. \hfill$\square$
\end{pf}
\begin{minipage}{\linewidth}
    \begin{defn}[Finite signed measures]
        Let $\mathcal{I}\subseteq\mathbb{R}$ be an interval and let $\mathcal{A}$ be a $\sigma$-algebra on $\mathcal{I}$. The space of finite signed measures on $(\mathcal{I},\mathcal{A})$ is defined as
        \begin{align*}
            \mathcal{M}(\mathcal{I}) &= \left\{\alpha:\mathcal{A}\to\mathbb{R}\;\middle|\;\alpha \text{ is countably additive, }|\alpha|<\infty\right\},
        \end{align*}
        where $|\alpha|$ denotes the total variation of $\alpha$ \cite[Chapters~1 and~6]{RudinRealAndComplex}.
    \end{defn}
\end{minipage}
Whenever the interval $\mathcal{I}$ on which the measure is defined is clear from context, we simply write $\mathcal{M}=\mathcal{M}(\mathcal{I})$.\\
For ease of notation we introduce the following operators that integrate the length-$\ell$ translates of a signal against a measure $\alpha \in \mathcal{M}([0,T-\ell])$:
\begin{itemize}
    \item Input operator:
    \begin{equation}
    \label{eqn:operatorinput}
        \begin{aligned}
            \mathcal{U}&:\mathcal{M}([0,   T-\ell])\to L^2([0, \ell];\mathbb{R}^{n_u})\\
            \alpha&\mapsto(\mathcal{U}\alpha)(t) =\int_0^{T-\ell}u_T(t+s)\mathrm{d}\alpha(s)
        \end{aligned}
    \end{equation}
    \item Output operator:
    \begin{equation}
    \label{eqn:operatoroutput}
        \begin{aligned}
            \mathcal{Y}&:\mathcal{M}([0,   T-\ell])\to L^2([0,\ell];\mathbb{R}^{n_y})\\
            \alpha&\mapsto(\mathcal{Y}\alpha)(t) =\int_0^{T-\ell}y_T(t+s)\mathrm{d}\alpha(s).
        \end{aligned}
    \end{equation}
    \item Initial state operator:
    \begin{equation}
        \begin{aligned}
            \mathcal{X}_0&:\mathcal{M}([0, T-\ell])\to \mathbb{R}^{n} \\
            \alpha&\mapsto\mathcal{X}_0\alpha =\int_0^{T-\ell}x_T(s)\mathrm{d}\alpha(s)
        \end{aligned}
    \end{equation}
    \item Joint input-output operator:
    \begin{equation}
    \label{eqn:IOOperator}
        \begin{aligned}
        \mathcal{H}&:\mathcal{M}([0,T-\ell])\to L^2([0,\ell];\mathbb{R}^{n_u+n_y})\\
        \alpha&\mapsto(\mathcal{H}\alpha)(t) =\int_0^{T-\ell}\begin{bmatrix}u_T(t+s)\\y_T(t+s) \end{bmatrix}\mathrm{d}\alpha(s)
    \end{aligned}
    \end{equation}
    \item Joint input-state operator:
    \begin{equation}
    \label{eqn:jointOperator}
        \begin{aligned}
            \mathcal{H}_x& : \mathcal{M}([0, T-\ell])\to L^2([0, \ell];\mathbb{R}^{n_u})\times\mathbb{R}^{n}\\
            \alpha&\mapsto(\mathcal{H}_x\alpha)(t) = \int_0^{T-\ell} \begin{bmatrix}u_T(t+s)\\x_T(s)\end{bmatrix} \mathrm{d}\alpha(s)
        \end{aligned}
    \end{equation}
\end{itemize}
The following theorem is the main result of this paper.

\begin{minipage}{\linewidth}
\begin{thm}
\label{thm:finitelengthbehaviors}
    Let $(u_T,y_T)\in\mathcal{B}_T$, where $u_T$ is $\ell$-PE of order $n$ and the system is controllable. Then, the closure of the range of the operator $\mathcal{H}$ in \eqref{eqn:IOOperator} equals the behavior of the system, i.e.,
    \begin{align}
        \overline{\operatorname{range}}(\mathcal{H}) &= \mathcal{B}_\ell.
    \end{align}
\end{thm}
\end{minipage}

The remainder of this section proves Theorem~\ref{thm:finitelengthbehaviors} in two steps:
\begin{itemize}
    \item[(i)] In Section~\ref{sec:rangeinB}, we show that every $\mathcal{H}\alpha$ is a trajectory of the system, i.e., $\overline{\operatorname{range}}(\mathcal{H}) \subseteq \mathcal{B}_\ell$.
    \item[(ii)] In Section~\ref{sec:Binrange}, we show that the range of $\mathcal{H}$ is dense in $\mathcal{B}_\ell$, i.e., every $(u_\ell, y_\ell) \in \mathcal{B}_\ell$ lies in the closure of $\operatorname{range}(\mathcal{H})$. We do this by showing that the adjoint of the joint input-state operator $\mathcal{H}_x$ has a trivial kernel.
\end{itemize}
Theorem~\ref{thm:finitelengthbehaviors} states that a pair of functions $(u_\ell, y_\ell)$ on $[0,\ell]$ is a trajectory of the system if and only if it can be approximated arbitrarily well by functions $\mathcal{H}\alpha$. The behavior is thus represented entirely by the measured data, which makes the theorem the continuous-time counterpart of the \textit{fundamental lemma} by \cite{WILLEMS2005325}. In discrete time, trajectories are generated by multiplying a Hankel matrix with a coefficient vector, i.e., by forming finite linear combinations of shifted data. In continuous time, the shift $s$ ranges over an interval, and $\mathcal{H}$ combines the shifted data by integrating against a measure $\alpha$. The integral operator $\mathcal{H}$ in \eqref{eqn:IOOperator} is the continuous-time analogue of the Hankel matrix and is known as Hankel integral operator. Finite signed measures include both weighted integrals of the translates (measures with a density) and finite linear combinations of individual translates (sums of Dirac measures).
\begin{rem}
    The operator in \eqref{eqn:IOOperator} is compact since $u_T\in  L^2([0,T];\mathbb{R}^{n_u})$ and $y_T\in L^2([0,T];\mathbb{R}^{n_y})$ ensure continuity of the map from each shift $s$ to the corresponding translated function in $L^2([0,\ell];\mathbb R^{n_u+n_y})$ \cite[Theorem~9.5]{RudinRealAndComplex}. A compact linear operator between Banach spaces cannot have a closed range of infinite dimension \cite[Theorem~4.18]{Rudin}. Therefore, the closure in Theorem~\ref{thm:finitelengthbehaviors} cannot be omitted.
\end{rem}
The excitation condition in Definition~\ref{def:PE} is weaker than existing notions in terms of regularity but stronger in terms of richness. The rank condition of \cite{9992347} and the notion of \cite{9882335} require (piecewise) differentiable data. However, they can be met by signals of finite complexity, such as sums of finitely many exponentials. Such signals are never $\ell$-PE, since their translates lie in a finite-dimensional subspace and therefore cannot be complete (Lemma~\ref{lem:denseinput}). This richness is not an artifact of Definition~\ref{def:PE}: completeness of the translates is necessary for the conclusion of Theorem~\ref{thm:finitelengthbehaviors}. Existing works avoid this requirement because they do not combine the shifted data with a fixed measure. In \cite{9992347}, the coefficients are time-varying and solve a differential equation. In \cite{RAPISARDA2023105603}, Chebyshev coefficients of the data are combined subject to derivative-consistency constraints.
\subsection{Step (i): Every $\alpha$ defines a system trajectory}
\label{sec:rangeinB}
Let us fix a measure $\alpha \in \mathcal{M}$ and let
\begin{align}
    (u_\ell, y_\ell) = \mathcal{H}\alpha.
\end{align}
We show that $(u_\ell, y_\ell) \in \mathcal{B}_\ell$, i.e., that $y_\ell$ is the output of \eqref{eqn:LTIsys} for the input $u_\ell$ and some initial state.\\
For $s \in [0,T-\ell]$ and $t \in [0,\ell]$, the state at time $t+s$ follows from the state at time $s$ as
\begin{align}
\label{eqn:stateevolve}
    x_T(t+s) &= e^{At}x_T(s) + \int_0^t e^{A(t-\sigma)}Bu_T(\sigma+s)\mathrm{d}\sigma
\end{align}
as derived in Appendix~\ref{app:StateEqn}. The output equation in \eqref{eqn:LTIsys} gives
\begin{align}
\label{eqn:outputevolve}
    &y_T(t+s) = Cx_T(t+s) + Du_T(t+s)\\
    &=\!Ce^{At}x_T(s) \!+\! \int_0^t \!Ce^{A(t-\sigma)}Bu_T(\sigma+s)\mathrm{d}\sigma \!+\! Du_T(t+s) \notag.
\end{align}
Integrating \eqref{eqn:outputevolve} against $\alpha$ yields
\begin{equation}
    \begin{aligned}
        &{y}_\ell(t)= \int_0^{T-\ell}y_T(t+s)\mathrm{d}\alpha(s)\\
        &= Ce^{At}\int_0^{T-\ell}x_T(s)\mathrm{d}\alpha(s) \\
        &\quad\quad\quad+\int_0^{T-\ell}\int_0^t Ce^{A(t-\sigma)}Bu_T(\sigma+s)\mathrm{d}\sigma\mathrm{d}\alpha(s)\\
        &\quad\quad\quad+ \int_0^{T-\ell} Du_T(t+s)\mathrm{d}\alpha(s).
    \end{aligned}
\end{equation}
Since $u_T$ is square integrable on a finite interval and $\alpha$ has finite total variation, the double integral converges absolutely. By Fubini's theorem \cite[Theorem~8.8]{RudinRealAndComplex}, we may therefore exchange the order of integration, which gives
\begin{align}
\label{eqn:afterIntSwitch}
    &{y}_\ell(t) = Ce^{At}\underbrace{\int_0^{T-\ell}x_T(s)\mathrm{d}\alpha(s)}_{x_\ell(0)}\notag\\
    &+ \int_0^t Ce^{A(t-\sigma)}B\underbrace{\int_0^{T-\ell}u_T(\sigma+s)\mathrm{d}\alpha(s)}_{{u}_\ell(\sigma)}\mathrm{d}\sigma \notag\\
    &+  D\underbrace{\int_0^{T-\ell}u_T(t+s)\mathrm{d}\alpha(s)}_{{u}_\ell(t)}.
\end{align}
Hence, $y_\ell$ is the output of \eqref{eqn:LTIsys} for the input $u_\ell$ and the initial state $x_0 = \mathcal{X}_0\alpha$, so $\mathcal{H}\alpha \in \mathcal{B}_\ell$ for every $\alpha \in \mathcal{M}$. Since the behavior $\mathcal{B}_\ell$ in \eqref{eqn:IOBehavior} is closed \citep[Chapter~2.7]{sontag} and contains $\operatorname{range}(\mathcal{H})$, we conclude
\begin{align}
\label{eqn:rangeSubsetB}
    \overline{\operatorname{range}}(\mathcal{H}) \subseteq \mathcal{B}_\ell.
\end{align}
\subsection{Step (ii): Every system trajectory lies in the closure of the range}
\label{sec:Binrange}
It remains to show that the closure of the range of $\mathcal{H}$ in \eqref{eqn:IOOperator} equals the entire behavior $\mathcal{B}_\ell$. Since a trajectory is determined by its input and initial state, it suffices to show that $\mathcal{H}_x$ has dense range.
\begin{lem}
\label{lem:InputToState}
If $\overline{\operatorname{range}}(\mathcal{H}_x)=L^2([0,\ell];\mathbb{R}^{n_u})\times\mathbb{R}^n$,
then $\overline{\operatorname{range}}(\mathcal{H})=\mathcal{B}_\ell$.
\end{lem}
\begin{pf}
Let
\begin{align}
    \Phi:(u_\ell,x_\ell(0))\mapsto(u_\ell,y_\ell)
\end{align}
map input and initial state to input and output. $\Phi$ is bounded and linear, and by
\eqref{eqn:afterIntSwitch}, $\mathcal{H}=\Phi\circ\mathcal{H}_x$.
By definition of $\mathcal{B}_\ell$, $\operatorname{range}(\Phi)=\mathcal{B}_\ell$.
Since $\Phi$ is continuous,
\begin{align}
    \mathcal{B}_\ell=\Phi\big(\overline{\operatorname{range}}(\mathcal{H}_x)\big)
    \subseteq\overline{\Phi(\operatorname{range}(\mathcal{H}_x))}
    =\overline{\operatorname{range}}(\mathcal{H}).
\end{align}
The reverse inclusion is \eqref{eqn:rangeSubsetB}.\hfill$\square$
\end{pf}
By Lemma~\ref{lem:InputToState}, the proof of Theorem~\ref{thm:finitelengthbehaviors} reduces to showing that
\begin{align}
    \overline{\operatorname{range}}(\mathcal{H}_x)&= L^2([0,\ell];\mathbb{R}^{n_u})\times\mathbb{R}^{n}.
\end{align}
To verify density, we use the adjoint operator $\mathcal{H}_x^\ast$. Since the codomain of $\mathcal{H}_x$ is a Hilbert space, the closure of the range is the orthogonal complement of the kernel of the adjoint \cite[Chapters~4~and~12]{Rudin},
\begin{align}
    \overline{\operatorname{range}}(\mathcal{H}_x) = (\operatorname{ker}(\mathcal{H}_x^\ast))^\perp.
\end{align}
If $\mathcal{H}_x^\ast$ has a trivial kernel, its orthogonal complement is the entire domain of $\mathcal{H}_x^\ast$ \cite[Chapter~3.2]{Neerven}, so the range of $\mathcal{H}_x$ is dense in its codomain. Together with Lemma~\ref{lem:InputToState}, we obtain
\begin{align}
    \operatorname{ker}(\mathcal{H}_x^\ast) = \{0\} \;  \implies\; \overline{\operatorname{range}}(\mathcal{H}) = \mathcal{B}_\ell.
\end{align}
\subsubsection{The adjoint operator}
We now compute the adjoint $\mathcal{H}_x^\ast$. By the Riesz representation theorem \cite[Theorem~4.12]{RudinRealAndComplex}, we identify the Hilbert space $L^2([0,\ell];\mathbb{R}^{n_u}) \times \mathbb{R}^n$ with its dual, using the inner product
\begin{align}
\label{eqn:innerProd}
    \langle (f,x), (v,\eta) \rangle = \int_0^\ell f(t)^\top v(t)\mathrm{d}t + x^\top\eta
\end{align}
for $f,v \in L^2([0,\ell];\mathbb{R}^{n_u})$ and $x,\eta \in\mathbb{R}^{n}$.
The adjoint is then the operator
\begin{align}
    \mathcal{H}_x^\ast : L^2([0,\ell];\mathbb{R}^{n_u}) \times \mathbb{R}^n
    \to \Big(\mathcal{M}([0,T-\ell])\Big)^\ast
\end{align}
defined by the duality pairing
\begin{align}
\label{eqn:canonicalpairing}
    \langle \mathcal{H}_x\alpha, (v,\eta) \rangle  = \langle \alpha, \mathcal{H}_x^\ast(v,\eta) \rangle_{\mathcal{M},\mathcal{M}^\ast},
\end{align}
for all $\alpha\in\mathcal{M}([0,T-\ell])$. Inserting \eqref{eqn:jointOperator} into the left-hand side of \eqref{eqn:canonicalpairing} and exchanging
the order of integration by Fubini's theorem gives
\begin{align}
    &\langle\mathcal{H}_x\alpha,(v,\eta)\rangle \notag\\
    &=\!\Big\langle\!\int_0^{T-\ell}\!u_T(\cdot\!+\!s)\mathrm{d}\alpha(s),v\!\Big\rangle\!+\!\Big(\int_0^{T-\ell}\!x_T(s)\mathrm{d}\alpha(s)\Big)^\top\!\eta \notag\\
    &=\!\int_0^{T-\ell}\!\Big(\!\int_0^\ell u_T(t+s)^\top v(t)\mathrm{d}t+x_T(s)^\top\eta\!\Big)\mathrm{d}\alpha(s). \label{eqn:deriveAdjoint}
\end{align}
The function in parentheses in \eqref{eqn:deriveAdjoint} is continuous and bounded on $[0,T-\ell]$, so integrating it against $\alpha$ defines a bounded linear functional on $\mathcal{M}$. Comparing with \eqref{eqn:canonicalpairing} identifies the adjoint as
\begin{align}
    \label{eqn:adjointOperator}
    (\mathcal{H}^\ast_x&(v,\eta))(s) = \int_0^\ell u_T(t+s)^\top v(t)\mathrm{d}t  + x_T(s)^\top\eta .
\end{align}
The first term is the cross-correlation of $u_T$ with $v$, which is the function $y_{v,0}$ in \eqref{eqn:integralPE}. The second term is the state trajectory $x_T$ weighted by $\eta$.

To show that the adjoint operator has a trivial kernel, we need to show that
\begin{align}
    \mathcal{H}_x^\ast(v,\eta) =0\quad \implies\quad \eta = 0 \quad &\text{and} \quad v =0.
\end{align}
The second term of the adjoint operator in \eqref{eqn:adjointOperator} is the
output $\tilde{y}$ of the \gls{lti} system
\begin{equation}
\label{eqn:virtualsys}
    \begin{aligned}
    \dot{x}_T(s) &= Ax_T(s) + B u_T(s), \quad x_T(0) \in \mathbb{R}^{n}\\
    \tilde{y}(s) &= \eta^\top x_T(s)
\end{aligned}
\end{equation}
for any vector $\eta\in\mathbb{R}^{n}$, which is an $n^{th}$-order strictly proper system driven by $u_T$. Since $u_T$ is $\ell$-\gls{pe} of order $n$, Definition~\ref{def:PE} with $q = 0$ implies that
\begin{align}
    y_{v,0}(s) &= \int_0^\ell u_T(t+s)^\top v(t)\mathrm{d}t
\end{align}
cannot equal such an output unless $v = 0$. We conclude that the ranges of the linear operators intersect trivially, i.e.,
\begin{align*}
    &\operatorname{range}\left(\eta\mapsto\eta^\top e^{As}x_T(0) + \eta^\top\int_0^s e^{A(s-\sigma)}Bu_T(\sigma)\mathrm{d}\sigma\right) \notag\\
    &\cap\quad\operatorname{range}\left(v(t)\mapsto\int_0^\ell u_T(t+s)^\top v(t)\mathrm{d}t\right) =\{0\}.
\end{align*}
If $\mathcal{H}_x^\ast(v,\eta) = 0$, then $y_{v,0} = (-\eta)^\top x_T$ lies in both ranges, so both terms vanish separately. In the following, we show that the first term vanishes only for $v = 0$ (Case~1) and the second only for $\eta = 0$ (Case~2).
\subsubsection*{Case 1: The cross-correlation term}
Suppose the first term vanishes, i.e.,
\begin{align}
\label{eqn:trivialKerv}
    y_{v,0}(s) = \int_0^\ell u_T(t+s)^\top v(t)\mathrm{d}t
    = \langle u_T(t+s), v(t) \rangle = 0
\end{align}
for all $s \in [0,T-\ell]$. Then $v$ is orthogonal to every length-$\ell$ translate of $u_T$. By Lemma~\ref{lem:denseinput}, the translates span a dense subspace of
$L^2([0,\ell];\mathbb{R}^{n_u})$, so $v$ is orthogonal to the whole space and hence $v = 0$.
\subsubsection*{Case 2: The state term}
Suppose the second term vanishes, i.e.,
\begin{align}
\label{eqn:etastate}
    \eta^\top x_T(s) = 0 \qquad \forall s \in [0,T-\ell].
\end{align}
To conclude that $\eta = 0$, we show that the state trajectory $x_T$ visits enough directions to span the whole state space, i.e., that $u_T$ excites all states of the system.
\begin{lem}
\label{lem:statespan}
    Let the pair $(A,B)$ be controllable and let $u_T$ be $\ell$-PE of order $n$. Then
    \begin{align}
        \operatorname{span}\{x_T(s) : s \in [0, T-\ell]\} = \mathbb{R}^{n}.
    \end{align}
\end{lem}
\begin{pf}
    Define
    \begin{align}
        V &= \operatorname{span}\{x_T(s):s\in[0,T-\ell]\}.
    \end{align}
    We show that $V$ contains the columns of $B$ and is invariant under $A$. Controllability then implies $V = \mathbb{R}^n$.\\
    Let $\mu \in V^\perp$. Then $\mu^\top x_T(s) = 0$ for all $s \in [0,T-\ell]$. Since $x_T$ is absolutely continuous, differentiating
    gives $\mu^\top\big(Ax_T(s) + Bu_T(s)\big) = 0$ for almost every $s$, i.e.,
    \begin{align}
    \label{eqn:qandetacondtion}
        q^\top u_T(s) = \eta^\top x_T(s)
    \end{align}
    with
    \begin{align}
        q = B^\top\mu \quad \text{and}\quad \eta = -A^\top\mu.
    \end{align}
    Suppose $q \neq 0$. The left-hand side of \eqref{eqn:qandetacondtion} is $y_{0,q}$ in \eqref{eqn:integralPE}, and the right-hand side is the output of the system
    \eqref{eqn:virtualsys} with output vector $\eta$. This contradicts Definition~\ref{def:PE} with $(v,q) = (0,q)$. Hence
    \begin{align}
        q &= B^\top\mu = 0
    \end{align}
    for every $\mu \in V^\perp$, i.e., every column of $B$ lies in $V$.\\
    With $q = 0$, \eqref{eqn:qandetacondtion} gives $\mu^\top A x_T(s) = 0$ . Hence $A^\top\mu \in V^\perp$ for every $\mu \in V^\perp$. Now let $x \in V$. For every $\mu \in V^\perp$, we obtain
    \begin{align}
        \mu^\top (Ax) = (A^\top\mu)^\top x = 0,
    \end{align}
    since $A^\top\mu \in V^\perp$ and $x \in V$. Hence $Ax$ is orthogonal to $V^\perp$, i.e., $Ax \in (V^\perp)^\perp = V$. Since $x \in V$ was arbitrary, $AV \subseteq V$.\\
    Since $V$ contains the columns of $B$ and is invariant under $A$, it contains the columns of $B, AB, \dots, A^{n-1}B$. By Assumption~\ref{ass:control}, these span $\mathbb{R}^n$, so $V = \mathbb{R}^n$.
    \hfill$\square$
\end{pf}
By Lemma~\ref{lem:statespan}, the vectors $x_T(s)$ span $\mathbb{R}^n$, so \eqref{eqn:etastate} implies $\eta = 0$. Together with Case~1, this shows that $\mathcal{H}_x^\ast$ has a trivial kernel, which completes the proof of Theorem~\ref{thm:finitelengthbehaviors}.

The proof uses Definition~\ref{def:PE} only for the true system and the true initial state. A weaker, system-dependent excitation condition is therefore sufficient for Theorem~\ref{thm:finitelengthbehaviors} (Appendix~\ref{app:systemPE}). Unlike Definition~\ref{def:PE}, it depends on the unknown state trajectory. Definition~\ref{def:PE} removes this dependence by requiring the condition for all $n^{th}$-order systems and all initial states.
\section{Illustrative example}
\label{sec:Example}
\begin{figure}
    \centering
    \includegraphics[width=\linewidth]{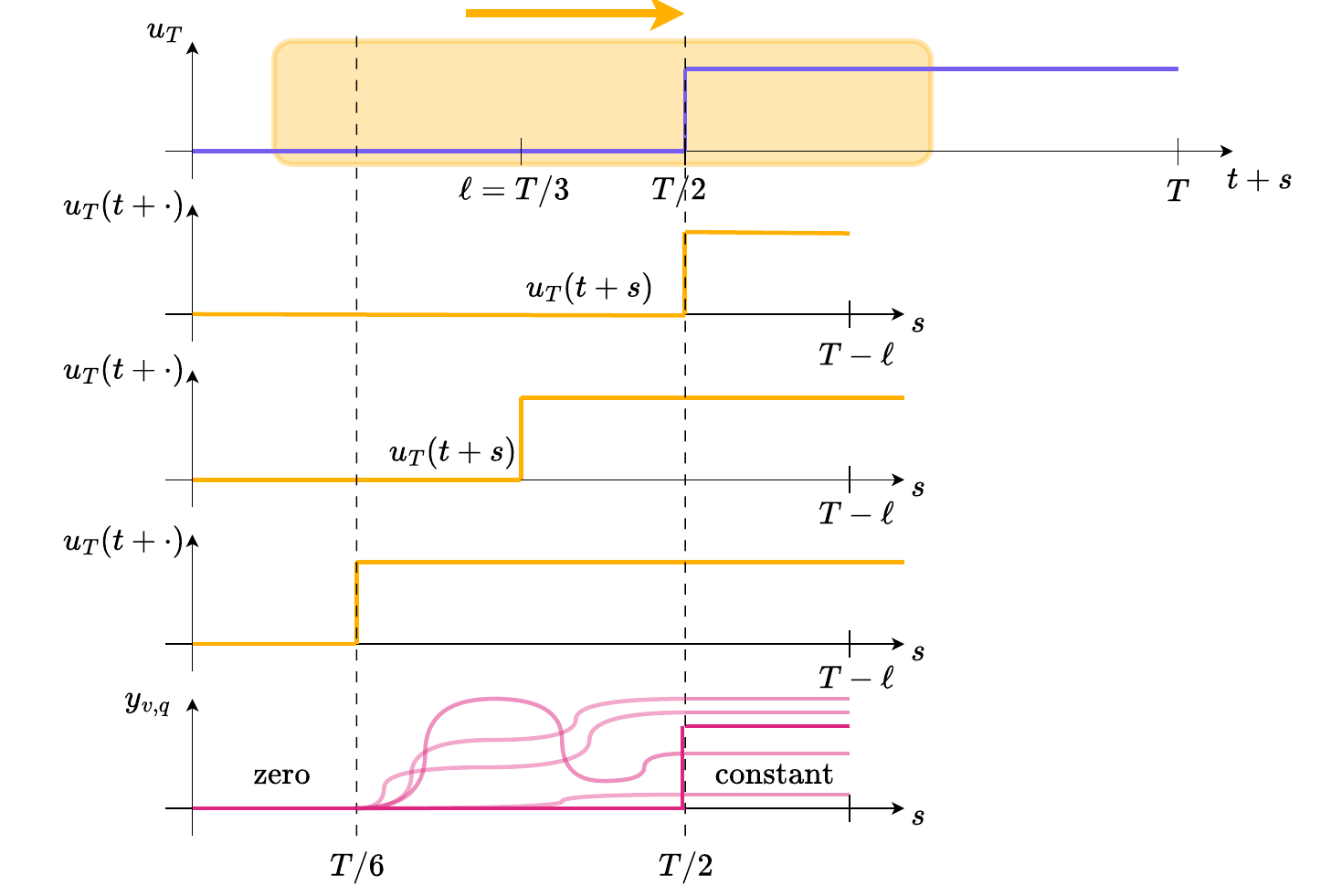}
    \caption{Calculation of $y_{v,q}$ (red) by integrating over shifted windows (yellow) of the input function $u_T$ (purple). The input signal is $\ell$-\gls{pe}.}
    \label{fig:PEInput}
\end{figure}

We illustrate the results for the first-order system
\begin{equation}
\label{eqn:minLTI}
    \begin{aligned}
        \dot{x}(t) &= ax(t) + bu(t)\\
        y(t) &= x(t)
    \end{aligned}
\end{equation}
with $a \neq 0$ and $b \neq 0$, so that the system is controllable. As input data, we apply the unit step at time $T/2$
\begin{align}
\label{eqn:StepInput}
    u_T(t) &= \begin{cases}0 \text{ for }t<\frac{T}{2}\\
    1 \text{ for }t\geq\frac{T}{2} 
    \end{cases}
\end{align}
on $[0,T]$ for some $T > 0$ as shown in purple in Figure~\ref{fig:PEInput}.
Initializing the system at $x_T(0) = 0$ results in the output
\begin{equation}
\label{eqn:StepOutput}
    \begin{aligned}
        y_T(t) &= \int_0^t e^{a(t-\sigma)}bu_T(\sigma)\mathrm{d}\sigma\\
         &= \begin{cases}0 \text{ for } t < \frac{T}{2}\\
        \frac{b}{a}(e^{a(t-T/2)}-1) \text{ for } t \geq \frac{T}{2}\end{cases}.
    \end{aligned}
\end{equation}
Throughout the example, we use the window length $\ell = T/3$, so that $T - \ell = 2T/3$ and $T/2 - \ell = T/6$.
\subsection{Persistency of excitation}
We first verify that the step input in \eqref{eqn:StepInput} is $\ell$-\gls{pe} according to Definition~\ref{def:PE}. For every fixed $t \in [0,\ell]$, the translate $s \mapsto u_T(t+s)$ is a unit step that jumps at $s = T/2 - t$, as shown by the yellow curves in Figure~\ref{fig:PEInput}. As $t$ increases from $0$ to $\ell$, the jump moves from $s = T/2$ to $s = T/2 - \ell = T/6$. The function $y_{v,q}$ in \eqref{eqn:integralPE} integrates these shifted steps weighted by $v(t)$ and adds $qu_T(s)$. This gives three regimes:
\begin{align}
\label{eqn:yvqStep}
    y_{v,q}(s) =
    \begin{cases}
        0, & s \in [0, T/6],\\[2pt]
        \int_{T/2-s}^{\ell} v(t)\mathrm{d}t, & s \in (T/6, T/2),\\[2pt]
        \int_0^{\ell} v(t)\mathrm{d}t + q, & s \in [T/2, T-\ell].
    \end{cases}
\end{align}
The function $y_{v,q}$ is shown in red in Figure~\ref{fig:PEInput}. In the first regime, no window has reached the step yet. In the second, the step lies inside the window, so $y_{v,q}$ depends on $v$. In the third, every window lies entirely after the step, so $y_{v,q}$ is constant.\\
Now suppose that $y_{v,q}$ is the output of a strictly proper LTI system of order $n$ driven by $u_T$. On $[0,T/6]$, both the input and $y_{v,q}$ are zero. Since the input remains zero until $T/2$, the output must remain zero on $[0,T/2)$. By \eqref{eqn:yvqStep}, this requires $v = 0$. Then $y_{0,q} = qu_T$ is zero and jumps to $q$ at $T/2$. This cannot be the output of any finite-order LTI system, whose output is continuous, unless $q = 0$. Hence, $u_T$ is $\ell$-\gls{pe} of every finite order.
\subsection{Generating a single trajectory}
We now use the data to generate a trajectory with the sinusoidal input $u_\ell(t) = \sin(\omega t)$ on $[0,\ell]$, for some $\omega > 0$. We choose $\alpha_0\in \mathcal{M}$ with density
\begin{align}
    \frac{\mathrm{d}\alpha_0(s)}{\mathrm{d}s} &= \omega\cos\Big(\omega\Big(\frac{T}{2}-s\Big)\Big)\mathbf{1}_{[\frac{T}{2}-\ell,\frac{T}{2}]}(s).
\end{align}
The input operator gives
\begin{align}
    (\mathcal{U}\alpha_0)(t) &= \int_{T/2-t}^{T/2} \omega\cos\Big(\omega\Big(\frac{T}{2}-s\Big)\Big)\mathrm{d}s
\end{align}
With the substitution $r = T/2 - s$, the boundaries transform as follows:
\begin{align}
    s=T/2-t\implies r=t,\quad s=T/2\implies r=0.
\end{align}
Consequently,
\begin{align}
\label{eqn:sinInput}
    (\mathcal{U}\alpha_0)(t) &= \int_0^t \omega\cos(\omega r)\mathrm{d}r = \sin(\omega t).
\end{align}
For the output operator, only shifts with $t + s > T/2$ contribute, since $y_T$ vanishes before the step. Substituting \eqref{eqn:StepOutput} in the output operator in \eqref{eqn:operatoroutput} gives
\begin{align}
    (\mathcal{Y}\alpha_0)(t) \!&=\! \frac{b\omega}{a}\!\int_{T/2-t}^{T/2}\!\big(e^{a(t+s-T/2)} \!-\! 1\big)\!\cos\Big(\!\omega\Big(\frac{T}{2}\!-\!s\Big)\!\Big)\!\mathrm{d}s \notag\\
    &= \frac{b\omega}{a}\int_0^t \big(e^{a(t-r)} - 1\big)\cos(\omega r)\mathrm{d}r.
\end{align}
Evaluating the integral yields
\begin{align}
    (\mathcal{Y}\alpha_0)(t)\!&= \!\frac{b\omega}{a}\left(\frac{a e^{at}\! -\! a\cos(\omega t) \!+\! \omega\sin(\omega t)}{a^2+\omega^2} \!-\! \frac{\sin(\omega t)}{\omega}\right) \notag\\
    &= \!\frac{b}{a^2+\omega^2}\!\big(\omega e^{at} \!-\! \omega\cos(\omega t) \!-\! a\sin(\omega t)\big).
    \label{eqn:partSolutio}
\end{align}
This is the zero-state response of \eqref{eqn:minLTI} to the input $\sin(\omega t)$. The initial state generated by $\alpha_0$ is $\mathcal{X}_0\alpha_0 = \int x_T(s)\mathrm{d}\alpha_0(s) = 0$, since $x_T = 0$ on the support of $\alpha_0$. Hence, $\mathcal{H}\alpha_0 \in \mathcal{B}_\ell$.
\subsection{Span of the behavior}
The measure $\alpha_0$ generates the input $\sin(\omega t)$  and results in the output with zero initial state. To obtain the same input with output signals starting from arbitrary initial states $x_0 \in \mathbb{R}$, we add a measure $\alpha_h \in \ker(\mathcal{U})$ that leaves the input unchanged. After the step at $T/2$, all input translates are identical. We therefore combine two translates after the step by constructing $\alpha_h$ as
\begin{align}
    \alpha_h(s) &= \frac{ax_0}{b(1-e^{aT/6})}\Big(\delta\Big(s-\frac{T}{2}\Big) - \delta\Big(s-\frac{2T}{3}\Big)\Big),
\end{align}
where $\delta$ denotes the Dirac measure. Since $u_T(t+T/2) = u_T(t+2T/3) = 1$ for all $t \in [0,\ell]$,
\begin{align}
    (\mathcal{U}\alpha_h)(t) &= \frac{ax_0}{b(1-e^{aT/6})}\big(u_T(t+\frac{T}{2}) - u_T(t+\frac{2T}{3})\big)\notag\\
    &= 0.
\end{align}
By linearity, the measure $\alpha = \alpha_0 + \alpha_h$ yields
\begin{align}
    (\mathcal U\alpha)(t)&= (\mathcal U(\alpha_0\!+\!\alpha_h))(t) = (\mathcal U\alpha_0)(t) + (\mathcal U\alpha_h)(t)\notag\\
    &= \sin(\omega t).
\end{align}
For the output, \eqref{eqn:StepOutput} gives
\begin{align}
    (\mathcal{Y}\alpha_h)(t) &= \frac{ax_0}{b(1-e^{aT/6})}\big(y_T(t+\frac{T}{2}) - y_T(t+\frac{2T}{3})\big) \notag\\
    &= \frac{ax_0}{b(1-e^{aT/6})}\frac{b}{a}e^{at}\big(1-e^{aT/6}\big) \notag\\
    &= x_0e^{at},\label{eqn:homoSolution} 
\end{align}
which is the free response of \eqref{eqn:minLTI} from the initial state $x_0$. By linearity and with \eqref{eqn:partSolutio} and \eqref{eqn:homoSolution}, the output is
\begin{equation}
    \begin{aligned}
        &(\mathcal Y\alpha)(t) \!=\! (\mathcal Y(\alpha_0\!+\!\alpha_h))(t) = (\mathcal Y\alpha_0)(t) + (\mathcal Y\alpha_h)(t)\\
        &= \! e^{at}x_0 \! + \! \frac{b}{a^2+\omega^2}\left(\omega e^{at} \! - \! \omega\cos(\omega t) \! - \! a\sin(\omega t) \right),
    \end{aligned}
\end{equation}
which is the response of \eqref{eqn:minLTI} to the input $\sin(\omega t)$ from the initial state $x_0$. Hence, every trajectory in $\mathcal{B}_\ell$ with input $\sin(\omega t)$ is exactly of the form $\mathcal{H}\alpha$.
\section{Conclusion}
\label{sec:Conclusion}
In this work, we introduced a notion of persistency of excitation for continuous-time \gls{lti} systems that is defined on the input data alone and does not require the input to be differentiable. The condition requires that no combination of the cross-correlation and the instantaneous value of the input can be reproduced by an $n^{th}$-order system driven by the same input. Under this condition, and for a controllable system, we showed that the closure of the range of a Hankel integral operator built from a single input-output trajectory equals the finite-horizon behavior of the system.\\
Future work will focus on leveraging the presented behavioral representation for data-driven control and system identification.
\begin{ack}                               
This work is part of Hollandse Kust Noord wind farm innovation program where CrossWind C.V., Shell, Eneco, Grow and Siemens Gamesa are teaming up; funding for the PhDs and PostDocs was provided by CrossWind C.V. and Siemens Gamesa.
\end{ack}

\bibliographystyle{plainnat}        
\bibliography{autosam}           
\appendix
\section{Derivation of the state equation}
\label{app:StateEqn}
For $s \in [0,T-\ell]$ and $t \in [0,\ell]$, the variation-of-constants formula with
initial time $s$ gives
\begin{align}
    x_T(t+s) = e^{At}x_T(s) + \int_s^{t+s} e^{A(t+s-\tau)}Bu_T(\tau)\mathrm{d}\tau.
\end{align}
Substituting $\tau = \sigma + s$ shifts the integration bounds from $[s,t+s]$ to
$[0,t]$ and yields \eqref{eqn:stateevolve}.
\section{A system-dependent excitation condition}
\label{app:systemPE}
The proof of Theorem~\ref{thm:finitelengthbehaviors} uses Definition~\ref{def:PE} only for the true system $(A,B)$ and the true initial state $x_T(0)$. The following weaker condition is therefore sufficient.
\begin{defn}
\label{def:systemPE}
    Let $T>\ell>0$ and $(u_T, x_T)\in\mathcal{B}_T^x$. The signal $u_T$ is said to be $\ell$-\gls{pe} with respect to $x_T$ if, for every pair $(v,q) \in  L^2([0,\ell];\mathbb{R}^{n_u})\times\mathbb{R}^{n_u}$ with $(v,q)\neq(0,0)$, the function $y_{v,q}$ in \eqref{eqn:integralPE} satisfies
    \begin{align}
        y_{v,q}(s) \neq \eta^\top x_T(s)
    \end{align}
    for every $\eta \in \mathbb{R}^n$ on a subset of $[0,T-\ell]$ of positive measure.
\end{defn}
Theorem~\ref{thm:finitelengthbehaviors} remains valid if Definition~\ref{def:PE} is replaced by Definition~\ref{def:systemPE}. However, Definition~\ref{def:systemPE} depends on the system matrices $A$ and $B$ and on the initial state $x_T(0)$ through the state trajectory $x_T$, which is not measured. It therefore cannot be verified from data when the system is unknown.
\end{document}